# Learning Outcomes supporting the integration of Ethical Reasoning into quantitative courses: Three tasks for use in three general contexts

Collaborative for Research on Outcomes and –Metrics;
Departments of Neurology; Biostatistics, Bioinformatics & Biomathematics; and Rehabilitation Medicine, Georgetown University, Washington, DC, USA

ORCID: 0000-0002-1121-2119

**Correspondence to:**
Rochelle E. Tractenberg
Georgetown University Neurology
Room 207, Building D
4000 Reservoir Rd., NW
Washington, DC, 20057 USA
**Email:** rochelle -dot- tractenberg -at -gmail -dot- com

Acknowledgement: There are no actual or potential conflicts of interest.

Running Head: Learning Outcomes (LOs) for Ethical Reasoning in Quantitative Courses

This synthesis was funded by a National Science Foundation Collaborative Grant (ER2 Standard Grant Award # 2220314)

This paper gives a brief overview of cognitive and educational sciences' perspectives on learning outcomes (LOs) to facilitate the integration of LOs specific to ethical reasoning into any mathematics or quantitative course. The target is undergraduate (adult) learners but these LOs can be adapted for earlier and later stages of learning. Core contents of Ethical Reasoning are: 1. its six constituent knowledge, skills, and abilities; 2. a stakeholder analysis; and 3. ethical practice standards or guidelines. These are briefly summarized. Five LOs are articulated at each of three levels of cognitive complexity (low/med/high), and a set of assignment features that can be adapted repeatedly over a term are given supporting these LOs. These features can support authentic development of the knowledge, skills, and abilities that are the target of ethical reasoning instruction in math and quantitative courses at the tertiary level. Three contexts by which these can be integrated are *Assumption* (what if the assumption fails?), *Approximation* (what if the approximation does not hold?), and *Application* (is the application appropriate? what if it is not?). One or more of the three core contents of Ethical Reasoning can be added to any problem already utilized in a course (or new ones) by asking learners to apply the core to the context. Engagement with ethical reasoning can prepare students to assume their responsibilities to promote and perpetuate the integrity of their profession across their careers using mathematics, statistics, data science, and other quantitative methods and technologies.

1. Getting "ethics" into quantitative courses: neither simple nor straightforward

Ethics is defined as "the principles of conduct governing an individual or a group (e.g., professional ethics)"[1]. “Ethics” is recommended content for statistics and data science curricula in higher education (e.g., American Statistical Association Undergraduate Guidelines Workgroup 2014; DeVeaux et al. 2017; National Academies 2018; Association of Computing Machinery Data Science Task Force, 2021), although notably not in undergraduate curricula in "business analytics" (Wilder & Ozgur, 2015) or mathematics (Saxe & Braddy, 2015), and also not in the European Union-focused EDISON Data Science Framework (Demchenko, Belloum & Wiktorski, 2017). One reason why some disciplinary curricular guidelines do not specify that "ethics" should be integrated might be that it is a very vague term. Even defined as professional guidelines, there may be more than one set that is relevant for an individual's practice (e.g., for federal statistics in the United States, a practitioner working with data might be guided by the *Data Ethics Tenets* (Office of Management and Budget, 2020), applicable to every federal employee working with data, or by the *Principles and Practices for Federal Statistical Agencies* (National Academies of Science, 2021) which are specifically utilized by designated agencies[2]. There are no alternative federal guidelines for mathematical practice in the United States.

In the context of data science statistical science, we can consider ethical practice standards that derive from statistics (American Statistical Association (ASA), 2022), computing (Association of Computing Machinery (ACM), 2018), and the area of specialization in which the data/statistical scientist is applying their statistical and computational expertise and methodologies. Critically, both the ASA and ACM assert the applicability of their ethical practice standards for *all* who utilize their domain knowledge, skills, and technologies. Whenever an individual -irrespective of membership in these professional organizations, degree or training, or job title- uses statistical practices or computing, the ASA and ACM ethical practice standards are relevant. “Upon entry into practice, all professionals assume at least a tacit responsibility for the quality and integrity of their own work and that of colleagues. They also take on a responsibility to the larger public for the standards of practice associated with the profession." (Golde & Walker, 2006: p. 10) This responsibility for the quality and integrity of data and statistical sciences work, is a responsibility to follow the professional ethics articulated for that type of work by these professional organizations. There is no single code of ethics for mathematical practices (but see American Mathematical Society, 2019 for a recently-updated version of their Code of Ethics; and Tractenberg et al., 2023 for *Mathematical Ethical Proto-Guidelines*).

---

[1] https://www.merriam-webster.com/dictionary/ethics

[2] The Principles and Practices apply specifically to 13 U.S. principal federal statistical agencies: Bureau of Economic Analysis (Department of Commerce); Bureau of Justice Statistics (Department of Justice); Bureau of Labor Statistics (Department of Labor); Bureau of Transportation Statistics (Department of Transportation); Census Bureau (Department of Commerce); Economic Research Service (Department of Agriculture); Energy Information Agency (Department of Energy); National Agricultural Statistics Service (Department of Agriculture); National Center for Education Statistics (Department of Education); National Center for Health Statistics (Department of Health and Human Services); National Center for Science and Engineering Statistics (National Science Foundation); Office of Research, Evaluation, and Statistics (Social Security Administration); and Statistics of Income (Department of Treasury). There are also three recognized federal statistical units: Microeconomic Surveys Unit (Federal Reserve Board); Center for Behavioral Health Statistics and Quality (Substance Abuse and Mental Health Services Administration; Department of Health and Human Services); and National Animal Health Monitoring System (Animal and Plant Health Inspection Service, Department of Agriculture).

It could be argued by instructors in any of the individual courses that make up undergraduate curricula in math, statistics, data science, and computing that theirs is not the most appropriate course into which "ethics" - especially not "ethical practice standards" - should be integrated. The argument might be that, since it is unclear if any of the students will in fact go on to be practitioners in the field, the professional practice standards are not yet relevant. Recognizing the potential for these objections, this paper hopes to empower instructors to overcome them. The solutions recognize that (and how) "ethics" is very difficult to teach -and even harder to assess in the typical undergraduate quantitative course. Quantitative-heavy courses do not facilitate engagement with "ethics" content. Moreover, instructors with quantitatively oriented materials would have to make room in the syllabus and in the course schedule to teach, and assess, engagement with ethics content. Within departments, each instructor might choose different teaching methods (e.g., discussion or case analysis), contents (e.g., case studies or online programs), and grading (e.g., rubrics; quizzes or tests). Thus, there are many problems that any instructor might encounter when considering the integration of ethics into a quantitative course or curriculum. These problems appear in Table 1. Solutions are also given in Table 1, and are elaborated throughout this paper.

Table 1. Problems with "integrating ethics" into quantitative courses, with solutions.

| **PROBLEM** | **SOLUTION(s)** |
|---|---|
| "Ethics" is recommended content for statistics and data science curricula in higher education (National Academies; ASA; NIH; NSF, etc.) | Teach ethical reasoning, not "ethics". Leverage – and teach – Ethical Practice Standards in codes and guidelines, not a topics list. Formulate & share specific learning outcomes. |
| "Ethics" is hard to teach! Harder to assess. | Formulate & share specific learning outcomes that promote ethical reasoning, not "ethics". Use stakeholder analysis, not "discussion". Leverage – and teach – Ethical Practice Standards, not a topics list. |
| Quantitative-heavy courses do not facilitate "ethics" content. | Use stakeholder analysis, not "discussion" based on homework problems and proofs. Leverage – and teach – Ethical Practice Standards, not a topics list. |
| Quantitative materials, and instructors, have to make room to teach, and assess, engagement with ethics content. | Leverage – and teach – Ethical Practice Standards, not a topics list. Maximize the case analysis/case study teaching method – to facilitate both teaching and assessing. |
| Each instructor, in different contexts, might choose different methods, contents, student work, grading. Creating consistency within/across programs is difficult. | Leverage – and teach – Ethical Practice Standards, not a topics list. Use stakeholder analysis, not "discussion". Teach ethical reasoning, not "ethics". |

In order to support the solutions listed in Table 1, some background is needed to contextualize the suggestions.

## 2. Background: Education Sciences

*What is a learning outcome (LO)*?
"Learning outcomes are statements of the knowledge, skills and abilities individual students should possess and can demonstrate upon completion of a learning experience or sequence of learning experiences." (Stanford University, n. d.)

https://web.stanford.edu/dept/pres-provost/irds/assessment/downloads/CLO.pdf.
**Figure 1.** Five phases of curriculum and instructional design (from Tractenberg et al. 2020, with permission; after Nichols, 2002).

LOs are a starting point for both instructional design (e.g., Nicholls, 2002; Diamond 2008; Nilson 2016; Tractenberg et al. 2020) and considerations of the integration of new content or ideas into existing curricula or courses.

*Role of LOs in curriculum and instructional design*
Five phases of curriculum and instructional development (Nicholls, 2002) are widely recognized:

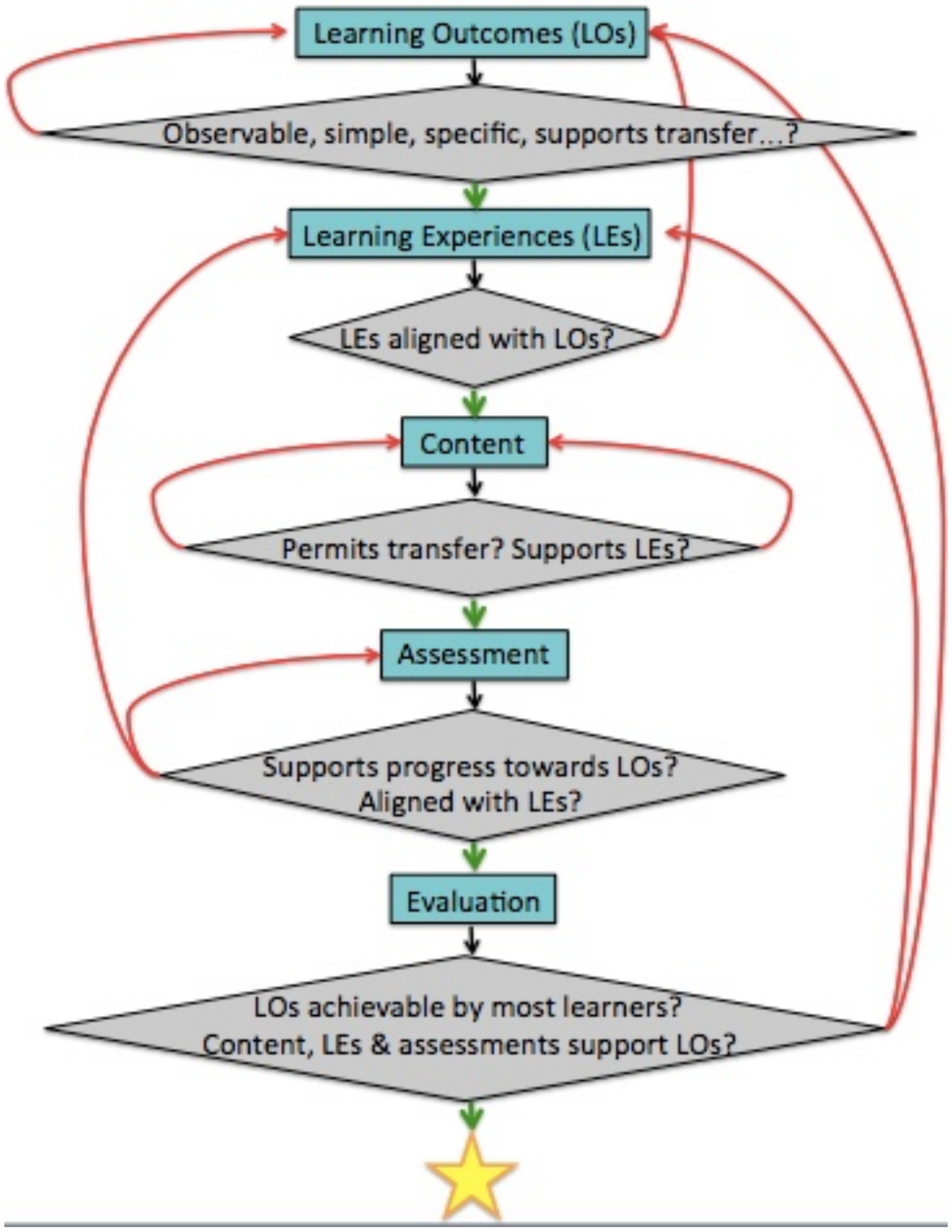

1. Identify aims and learning outcomes (LOs);
2. Identify/create learning experiences (LEs, lectures, exercises, readings, etc.) that will help students achieve the aims and outcomes;
3. Select content that is relevant to outcomes;
4. Identify/develop assessments to ensure learner is progressing towards outcomes;

5. Evaluate the effectiveness of the learning experiences for leading/developing learners to the outcomes.

This figure, showing the five phases of curriculum and instructional design (in blue boxes), also shows how these phases are iteratively informed by LOs. In fact, LOs drive all decisions in curriculum and in course development (Tractenberg et al. 2020; Nilson 2016).

Nilson (2016) noted, "your student learning outcomes provide the foundation for every aspect of your course, and you should align all the other components with them." (p. 129). Note that even course content is derived indirectly from LOs as a function of the learning experiences (teaching and classroom activities). Figure 1 shows how, at each phase, support for LOs is key in decision making about all other aspects of instruction (including assessment of learning and evaluation of teaching).

*Formulating LOs*

The following serves as a template for formulating a learning goal/outcome: fill in the blank with a specific, and concretely observable verb phrase that is relevant to the course/program: “At the end of the course, students should be able to ___.” Bloom's taxonomy of cognitive behaviors can be used to ensure that the LOs represent observable, evaluable (i.e., gradable) targets for learners as well as for instructors. Bloom's taxonomy outlines six levels of cognitive complexity, published in 1956 (elaborated in Tractenberg et al. 2013). Considerations about how learners can demonstrate (perform) at each Bloom's level appear in the list below.

1. *Remember/Reiterate*- performance is based on recognition of previously seen example
2. *Understand/Summarize*- performance summarizes info already known/given
3. *Apply/Illustrate* - performance extrapolates from seen examples to (really) new examples by applying rules - correctly
4. *Analyze/Predict* - performance requires analysis and prediction, using rules - judiciously
5. *Create/ Synthesize* –performance yields something innovative and novel; creating, describing, and justifying something new.
6. *Evaluate/Compare/Judge* - performance involves the application of guidelines, rather than rules, and can involve subtle differences arising from comparison or evaluation of abstract, theoretical, or otherwise not-rule-based decisions, ideas, or materials.

Cognitive complexity is reflected in many different taxonomies (see e.g., Moseley et al. 2005). Bloom's taxonony is perhaps more familiar to US or North American instructors than to those in other parts of the world. A superficial, but still important, utility of Bloom's taxonomy is the myriad lists of verbs that have been created and publicized to facilitate the creation of LOs. Filling in the blank in the template LO formulation above, with a verb that is appropriate to the course as well as the student's abilities, will generate evaluable - and achievable - LOs. The importance of cognitive complexity organized hierarchically like Bloom's taxonomy is also important in formulating LOs because:

a. It is important to recognize that when LOs have high complexity as their focus, they also assume, and require, engagement with the material at all the *earlier/lower* levels of complexity on the way to the target high level. The complexity level of any LO must be consistent with both the learners' abilities and the instruction and assessment that are planned to support and determine whether the LO is acheived. It is unlikely that an LO articulated a higher complexity level compared to the instruction can be achieved.

b. The same (well-formed) LO can be repurposed throughout a course with Bloom's level verbs increasing steadily, so that learners begin with low-cognitive complexity versions/activities and progress to the target (higher) level in a coherent way. This is important for ensuring that the learner recognizes how their own performance is changing, enabling them to monitor their learning and become more self-sufficient at metacognitive activities (like self-directed learning) (Tractenberg et al. 2017).

*Scaffolding and Blooms*

Scaffolding refers to both the gradual increase in independence of learners as they engage with activities over time (reducing the support from the instructor; Wells & Edwards, 2015; see also Vygotsky, 1978) and also to building a developmental trajectory for learners to advance in their sophistication (Tractenberg, 2017). This can be efficiently accomplished by re-using a task to revisit/reinforce or build towards achieveing LOs with increasing cognitive complexity over time. Additional LOs based on the same task/material efficiently build towards new and more sophisticated targets by increasing cognitive complexity of behavior students must exhibit. Table 2 shows one example of the multitudes of lists of concrete and observable behaviors students can be asked to perform at each level.

| Knowledge | Comprehension | Application | Analysis | Synthesis | Evaluation |
|---|---|---|---|---|---|
| duplicate | classify | apply | break down | assemble | argue |
| identify | describe | construct | classify | categorize | conclude |
| know | discuss | dramatize | differentiate | compose | criticize |
| list | explain | interpret | compare | design | defend |
| match | give examples | practice | contrast | modify | estimate |
| memorize | paraphrase | produce | distinguish | reconstruct | justify |
| recite | restate | solve | outline | revise | predict |
| repeat | reword | use | separate | summarize | support |

Melissa A Nelson 2011

**Table 2**. Bloom's original (1956) 6-level taxonomy of cognitive complexity (column headings) with eight example verbs appropriate at each. Downloaded from Melissa A. Nelson, 2011 https://manelsonportfolio.blogspot.com/search/label/Bloom

Table 2 supports the claim that LOs featuring - requiring- high level Bloom's level cognitive complexity will not be achievable by those who are taught to, or assessed in their ability to, perform at lower Bloom's level complexity. In considering how to scaffold around new LOs, and how to integrate ethical reasoning into quantitative courses in particular in a realistic and achievable way, instructors can leverage Bloom's taxonomy to ensure that the reduction in support that is characteristic of instructional scaffolding promotes achievement of the target LOs.

*Achieving psychometric validity through instructional design*

Messick (1994) outlines three important questions for evaluating the *effects* of instruction (i.e., whether learning goals were met):

1. What are the knowledge, skills, and abilities (KSAs) the instruction should lead to?
2. What actions/behaviours by the students will reveal these KSAs?
3. What tasks will elicit these specific actions or behaviours (and so reveal whether students did learn the target KSAs)?

NB: **Substitute "LO" for "KSA" in Messick's questions, because the effects of instruction are the achievement of learning outcomes (ostensibly)**.

*How to use Messick's three questions:*

Note that assessment (determination of whether LOs were met) is a core element of curriculum and instructional development (Hutchings 2016; National Institute for Learning Outcomes Assessment, 2016; Tractenberg, 2021). Messick, a psychometrician, focused on ensuring that assessment is valid - i.e., supports claims that "material was learned". Valid assessment requires the LOs to be articulated, in addition to specifying what students should be able to do to show that those LOs were achieved. This is one reason why Bloom's verbs - all concrete and observable - are so critical in LO formulation. Bloom's helps you devise LOs that are specific to what the learner actually can do with the new information/abilities you're going to teach them. That helps you answer Messick's first question (Messick #1). Bloom's verbs can make Messick #2 explicit, and also support the creation of tasks or assessments that accomplish Messick #3.

Once you have an outline of the scaffold for how to get learners from low to high Bloom's complexity, you can plan your time to teach towards those LOs, and give sufficient opportunities for learners to practice them with feedback, and then to demonstrate the level you targeted. So, while it seems a lot of background, taking Bloom's complexity and Messick's assessment features into account can make scaffolding more straightforward and achievable. Ultimately this can make instructional development more efficient; considerations of the elements that assessments should have for optimal utility can support informative and useful assessments (Tractenberg, 2021).

## 2. Contexts in quantitative courses for three activities or tasks that comprise Ethical Reasoning

There are three activities or tasks that can be utilized to integrate instruction about ethical reasoning into any math course: working with Guidelines (GLs); Stakeholder Analysis; and at least the first of the six ethical reasoning KSAs. For quantitative courses, there are three general contexts into which ethical reasoning and content can be introduced without significant deviation from the core content. These contexts are **assumptions, approximations,** and **applications**. Specifically, many mathematical and quantitative courses involve at least some discussion of assumptions that must be met for proofs or other methods to be relevant; approximations that are introduced when constructs or methods are defined, and applications of methods or constructs to "real world" situations. Whenever any of these appears in a course, it represents an opportunity to integrate one or more of the ethical reasoning activities: ask students to consider what might happen if the assumptions are not met, the approximations do not hold, or the application is not appropriate. The next section discusses the three tasks/activities: 1) engage with Stakeholder Analysis; 2) work with/learn about Guidelines; and 3) learn and demonstrate Ethical Reasoning.

**Stakeholder Analysis**: *a template for considering the impact of decisions in mathematical and quantitative practice.*

A stakeholder is defined “one who is involved in or affected by a course of action” (Merriam-Webster https://www.merriam-webster.com/dictionary/stakeholder); "in the context of ethical case analysis, the stakeholder is simply an individual, or group, that might be affected by the outcome of the case." (Tractenberg, 2019). The “stakeholder analysis template” (Tractenberg, 2019) is a method for understanding the harms and benefits, as well as the stakeholders, who may be impacted whenever decisions are made. Originally conceptualized for instruction about ethical considerations in statistics and data science (see Tractenberg 2022-A, 2022-B), the Stakeholder Analysis was utilized in the series of each of a series of tasks common to the collection, analysis/manipulation, and drawing of inferences or conclusions based on data in any shape or size. The Stakeholder Analysis is in the process of being integrated into instruction for undergraduate mathematics courses supporting science, technology, engineering, and mathematics (STEM) disciplines (e.g., Tractenberg & Thornton, 2024).

| **Potential result:**<br><br>**Stakeholder[1]:** | **HARM[4,5]** | **BENEFIT[4,5]** |
|---|---|---|
| **YOU[2,3]** | | |
| **Your boss/client** | | |
| **Unknown individuals[2]** | | |
| **Employer** | | |
| **Colleagues** | | |
| **Profession** | | |
| **Public/public trust** | | |

**Table 3.** Stakeholder Analysis activity: Use this template to consider (list) harms and benefits – if assumptions don’t hold; approximations are incorrect; or applications are inappropriate. Table adapted from Tractenberg (2019), with permission.

Notes on Stakeholder Analysis table:

1. Knowing to whom harms may accrue can point to where the Guidelines can assist in decision making.
2. Recognize/match whether or not benefits or harms accrue to any stakeholders. Are some harms “worse” or some benefits “better”?
3. If there are no recognizable harms, and plausibly no “unknowable” harms <for which your decision would be responsible>, then there can be no conflict. Your failure to recognize something doesn’t mean it does not exist.
4. Learning how to use this table and complete a case analysis is essential for enabling informed decisions about ethical challenges.
5. All harms are not the same; all the benefits are not the same; and harms and benefits are not exchangeable.

The Stakeholder Analysis itself supports articulation and ultimately, achievement and documentation of achievement, of three learning objectives:

1. Describe how different individuals ("stakeholders") may be affected by decisions and actions;
2. Enumerate harms and benefits that are most clearly relevant for each stakeholder with respect to the activity; and
3. Identify which ethical guideline or policy principles, practices, and/or specific elements seem most relevant to the specific task.

The Stakeholder Analysis template is intended to facilitate teaching and learning of quantitative and computational ethical practice standards like the American Statistical Association (ASA) Ethical Guidelines for Statistical Practice (ASA, 2022) and the Association for Computing Machinery (ACM) Code of Ethics (ACM, 2018). It can also be used with other guideline documents such as the Ethical Guidelines of the American Mathematical Society (AMS 2019) or the Mathematics Ethical proto-Guidelines (Tractenberg et al. 2024)

***Ethical Guidelines:***

As noted, ethical practice standards exist (and have been maintained for decades) for statistics (ASA, 2022) and computing (ACM, 2018). For the ethical practice of mathematics, there are (as of 2024) no universally ackhowledged ethical practice standards, but there are multiple sources from which an instructor, student, or practitioner can glean actionable guidance. These are thematically analyzed in the following tables.

The **ACM Code of Ethics** (ACM 2018), has four core areas, with 2-9 elements in each (total = 24 narrative elements). The full Code is included in the Appendix; an outline of the Code follows the four sections:

1. General Moral Principles (7)
2. Professional Responsibilities (9)
3. Professional Leadership Principles (7)
4. Compliance with the Code (2)

The **ASA Ethical Guidelines for Statistical Practice** have a total of 72 elements organized under eight Principles and one Appendix for institutions and organizations (12 elements). The full Guidelines is included in the Appendix; an outline is given in Table 4.

**Table 4.** Principles of the ASA Ethical Guidelines for Statistical Practice (ASA, 2022)

| |
|---|
| A. Professional Integrity & Accountability: Professional integrity and accountability require taking responsibility for one's work. Ethical statistical practice supports valid and prudent decision making with appropriate methodology. The ethical statistical practitioner represents their capabilities and activities honestly, and treats others with respect. (12 elements) |
| B. Integrity of data and methods: The ethical statistical practitioner seeks to understand and mitigate known or suspected limitations, defects, or biases in the data or methods and communicates potential impacts on the interpretation, conclusions, recommendations, decisions, or other results of statistical practices. (7 elements) |
| C. Responsibilities to Stakeholders: Those who fund, contribute to, use, or are affected by statistical practices are considered stakeholders. The ethical statistical practitioner respects the interests of stakeholders while practicing in compliance with these Guidelines. (8 elements) |
| D. Responsibilities to research subjects, data subjects, or those directly affected by statistical practices: The ethical statistical practitioner does not misuse or condone the misuse of data. They protect and respect the rights and interests of human and animal subjects. These responsibilities extend to those who will be directly affected by statistical practices. (11 elements) |
| E. Responsibilities to members of multidisciplinary teams: Statistical practice is often conducted in teams made up of professionals with different professional standards. The statistical practitioner must know how to work ethically in this environment. (4 elements) |
| F. Responsibilities to Fellow Statistical Practitioners and the Profession: Statistical practices occur in a wide range of contexts. Irrespective of job title and training, those who practice statistics have a responsibility to treat statistical practitioners, and the profession, with respect. Responsibilities to other practitioners and the profession include honest communication and engagement that can strengthen the work of others and the profession. (5 elements) |
| G. Responsibilities of Leaders, Supervisors, and Mentors in Statistical Practice: Statistical practitioners leading, supervising, and/or mentoring people in statistical practice have specific obligations to follow and promote these Ethical Guidelines. Their support for – and insistence on – ethical statistical practice are essential for the integrity of the practice and profession of statistics as well as the practitioners themselves. (5 elements) |
| H. Responsibilities regarding potential misconduct: The ethical statistical practitioner understands that questions may arise concerning potential misconduct related to statistical, scientific, or professional practice. At times, a practitioner may accuse someone of misconduct, or be accused by others. At other times, a practitioner may be involved in the investigation of others' behavior. Allegations of misconduct may arise within different institutions with different standards and potentially different outcomes. The elements that follow relate specifically to allegations of statistical, scientific, and professional misconduct. (8 elements) |
| APPENDIX: Responsibilities of organizations/institutions: Whenever organizations and institutions design the collection of, summarize, process, analyze, interpret, or present, data; or develop and/or deploy models or algorithms, they have responsibilities to use statistical practice in ways that are consistent with these Guidelines, as well as promote ethical statistical practice. (organizations 7 elements; leaders 5 elements; 12 elements total) |

Although they have diverse origins and organizing principles, the ASA and ACM ethical practice standards are highly concordant (Tractenberg, 2022-A).

The AMS Code of Ethics was thematically analyzed to identify a core set of ethical obligations for mathematical practice from the ASA (2018 edition) and ACM (2018) practice standards (Buell et al. 2022; Tractenberg et al. 2024).

**Table 5. AMS Code of Ethics thematic elements (Buell et al. 2022; used with permission):**

| **I. Mathematical research and its presentation** |
|---|
| Do not plagiarize, correct attribution when appropriate is essential. |
| Be knowledgeable in your field |
| Give appropriate credit |
| Do not claim a result in advance of its having been achieved; publish full details of results without unreasonable delay after announcing results. |
| A claim of independence may not be based on ignorance of widely disseminated results |
| Ensure appropriate authorship |
| **II. Social responsibility of Mathematicians** |
| Encourage and promote mathematical ability without bias and review programs to ensure consideration of a full range of students. |
| Avoid conflicts of interest and bias in reviewing, refereeing, or funding decisions. |
| Respect referee anonymity |
| Resist excessive secrecy, promote dissemination/publication |
| Disclose implications of work to employers and the public when work may affect public health, safety, or general welfare. |
| Do not exploit workers with temporary employment at low pay/excessive work) |
| **II. Education and Granting of Degrees** |
| Granting a degree means certifying competence for work. |
| PhD level work is ensured by the degree grantors to be high level and original |
| PhD is only awarded to those with sufficient knowledge outside the thesis area. |
| Degree grantors must honestly inform degree earners about job market/ employment prospects. |
| **IV. Publications** |
| Editors should be reasonably sure of the correctness of articles they accept. |
| Editors should ensure timely and current reviews. |
| Submissions for review are treated as privileged information. |
| Editors must prioritize the first submitted version of a paper. |
| Editors must inform authors if there is a delay in potential publication |
| Publication cannot be delayed for any reason except the authors' interest/actions. |
| Date of submission and revisions must be published with any article. |
| Editors must be given/accept full responsibility for their journals, resist outside agency pressures and notify the public of such pressure. |
| Editors and referees must respect the confidentiality of all submitted materials as appropriate. |
| Mathematical publishers must respect the mathematical community and disseminate work accordingly. |
| The American Mathematical Society will not play a role/endorse any research journal where any acceptance criterion conflicts with the principles of the AMS guidelines. |

**Mathematics Ethical proto-Guidelines** (Tractenberg et al. 2024)
After a thorough analysis of the guidance from the ASA (2018), ACM (2018), and AMS (2019) codes and guidelines, and the Code of the Mathematical Association of America, Tractenberg et al. 2024 derived a set of 52 elements for consideration by the mathematics community to constitute "Ethical Guidelines for the Practice of Mathematics". Of these 52 elements, only one was not endorsed as relevant for ethical mathematical practice by a majority (>50%) of the national sample of 142 respondents. Many of these items were derived from the lists above (see Buell et al. 2022). The final proto-Guidelines comprises 44 items, organized according to whether they describe responsibilities of the individual practitioner or the practitioner in a leader/mentor/supervisor/instructor role. For the individual practitioner, there are 12 General, 10 Profession, and 11 Scholarship items. For the practitioner in a leader/mentor/supervisor/instructor role, there are an additional 11 items (4 General/7 Professional). The full proto Guidelines for Ethical Mathematical Practice appears in the Appendix. These proto-Guidelines were derived from existing guidance for mathematics together with considerations of statistics, data science, and computation; while Buell et al. 2022 and Tractenberg et al. 2024 report the same survey of the community of mathematics practitioners in the US relating to these items, they are not finalized or generally acknowledged (as of October 2023). Instructors could utilize the proto-Guidelines, or they can use the results (also reported in Buell et al. 2022 and Tractenberg et al. 2024) of a thematic analysis of respondent suggestions for additions to the proto-Guidelines:

a. Workplace (not teaching, even if you do teach at work) – basic human respect/rights non-violations
b. Educating (if this is your primary job/role or if you teach only as part of mentoring/collaborating): teaching effectively, grading objectively; doing your best to promote learning.
c. Scholarship (writing, reviewing, and correcting errors, even if you made them) – respect for other's work and others' input to your work.
d. Respect for the profession/stewardship (apart from scholarship, ethical & objective reviewing)
e. Math in the world: effectively preparing learners & users of math; not gatekeeping.
f. Recognizing and effectively/respectfully treating stakeholders in work, teaching, scholarship, use of math, and the profession.

If instructors do not want to integrate as-yet-not-universally-acknowledged proto-Guidelines into their courses, these suggested themes could be useful for engaging students because they suggest an organizational scheme for "what constitutes ethical mathematical practice". To integrate ethical considerations into a mathematics course, these themes could be utilized in conjunction with the Stakeholder Analysis template - for example, considering whether harms or benefits accrue to diverse stakeholders when each of these themes is ignored. If scholarship is not carried out respectfully, do harms accrue to the employer (workplace) or to the public (math in the world)? It could be argued that, if an error is detected, but not publicly corrected, then a benefit arising from not being respectful in scholarship *could* accrue to oneself and an employer (benefit = no one knows you/the company made an error), while this same behavior creates potential harm to the public if anyone uses - or builds new work on - the incorrect result. A similar harm accrues to the profession and wider scientific community; these harms and benefits would appear in the Stakeholder Analysis as each of the themes is considered in a task like this. Note that such a task could be made to fit a lower Bloom's level if it was a matching or fill in the blank task; or, a high Bloom's task if it was an essay prompt.

***Ethical Reasoning****: six knowledge, skills, and abilities (KSAs), plus stakeholders, and guidelines.*

Ethical Reasoning is a 6-step process (Tractenberg & FitzGerald, 2012) that also involves understanding what our ethical obligations are, and recognizing stakeholders and harms and benefits that may accrue to them in our practice (in step 1) (Tractenberg 2022-A, 2022-B). Ethical challenges can be recognized whenever behavior is not consistent with law, policy, or ethical guidelines. However, ethical challenges may also arise when some stakeholders are harmed while others benefit. Therefore, ethical guidelines and stakeholder analysis comprise core foundational information. By articulating the potential impacts on stakeholders of our decisions (or contexts where assumptions and approximations fail to hold or applications are inapprorpriate or unjustified), and by ensuring that our fundamental quantitative practices are consistent with the community norms such as are as described in the Ethical Guidelines for Statistical Practice (ASA, 2022), ACM Code of Ethics (ACM, 2018), Ethical Guidelines of the American Mathematical Society (2019) or Mathematics Ethical proto-Guidelines (Tractenberg et al. 2024), we assure ourselves and others that our quantitative practice is ethical. However, whenever an ethical challenge arises, it must be *recognized*, and a decision must be made.

To *make* and then support the decision, you must:

1. Identify/ assess your prerequisite knowledge – using codes, guidelines, or other policy/law as well as a stakeholder analysis
2. Identify relevant decision-making frameworks (e.g., virtue or utilitarianism, although these are not the only frameworks)
3. Recognize an ethical issue (<u>and that a decision must be made</u>) –codes or guidelines will help to pinpoint

4. Identify and evaluate alternative actions -codes or guidelines might specify at least one course of action
5. Make & justify a decision
6. Reflect on the decision

These are the six steps in ethical reasoning, and the six elements of knowledge, skill, and ability that are outlined in the Mastery Rubric for Ethical Reasoning (Tractenberg & FitzGerald 2012).

Note. If there is no ethical challenge - i.e., the individual just wants to make sure they practice ethically, then Ethical Reasoning is limited to KSA 1, which is still quite a lot of information (Tractenberg 2022-A). You can use ethical guidelines or codes to ensure that you just do your work ethically, but in that case, you are not making a decision about what to do in the face of an ethical challenge (which you only do when there IS ONE, i.e., you or someone else has done something about which you need to decide what to do). If you complete KSAs 1-3 for a situation and recognize that, while it might require a decision to be made by you, it does not violate or appear to violate any ethical guideline, then the actions you can take (for KSA 4) are things like, "note this situation for continued monitoring" or "note this for future consideration".

In summary, there are three tasks that can be utilized to integrate instruction about ethical reasoning into any math course: working with codes and guidelines; Stakeholder Analysis; and at least the first of the six ethical reasoning KSAs (which combines the codes/guidelines and Stakeholder Analysis). These tasks can be executed in the context of assumption failures, approximation failures, and unjustified or inappropriate application of methods. This generates a 3x3 matrix of opportunities to elicit evidence that LOs are achieved (Messick #2). All that remains is to leverage what we know about cognition, plus course specifics, to construct appropriate LOs. Note that the integration of ethics content can be meaningful and directed at increasing responsible and ethical practice in statistics, data science, mathematics, and computing *without case analysis*. (Case analysis would involve all six Ethical Reasoning KSAs, in essay form, in response to a given vignette (e.g., Tractenberg 2022-B).

## 3. Learning outcomes featuring Stakeholder Analysis, codes/guidelines, and ethical reasoning at three Bloom's levels of complexity

### *Learning outcomes early in course (low Blooms)*

Overall teaching objective: orientation to ethical by introducing “prerequisite knowledge": what is a stakeholder, what is the relevant code/guidelines, and what is ethical reasoning. NB: for those programs or instructors who intend to teach all of the ethical reasoning KSAs, and especially where the Mastery Rubric for Ethical Reasoning should be shared with learners so they can see where the instruction is (intending) to lead them.

The *teaching* objectives are to promote the creation of a new schema for the student – that ethical reasoning is a learnable, improvable set of knowledge, skills, and abilities; and that they can develop, and self-assess this development, of the KSAs of ethical reasoning in progressive – manageable- ways. Examples of learning outcomes fitting Stanford LO-writing criteria, include:

After this section of the course, students will:

1. Map a vignette into the pre-analysis elements of a case analysis utilizing the relevant Ethical Reasoning KSAs. (*laying the foundation for case analysis, introducing* Ethical Reasoning *KSAs*)
2. Identify elements from the specific code/guideline that are relevant whenever assumptions, approximations, or preconditions are not met or mathematical applications are not justified/appropriate across a variety of mathematical contexts. (*establishing familiarity with the target code/guideline and how they function as a normative part of quantitative practice)*
3. Identify harms and benefits to stakeholders whenever assumptions, approximations, or preconditions are not met or mathematical applications are not justified/appropriate across a variety of mathematical contexts. (*introducing Stakeholder Analysis as a normative part of math*)
4. Summarize the relationships between stakeholder harms and benefits and the code/guidelines. Outline how stakeholder analysis and codes/guidelines act as prerequisite knowledge in ethical reasoning. (*reinforcing both* Ethical Reasoning *KSAs and engagement with how the* Ethical Reasoning *steps work with the code/guidelines and the Stakeholder Analysis*)
5. Describe how completing the tasks (below) represent evidence of their performance of each Ethical Reasoning KSA, including KSAs they have not learned/demonstrated yet. (*reinforcing engagement with* Ethical Reasoning *KSAs and what evidence of their learning looks like*)

The tables outline three tasks across three contexts specific to undergraduate mathematics courses, with structural changes reflecting increased cognitive complexity over iterations

**Table 6.** Formula for creating an assignment (task) where learners must perform at Bloom's level 1-2: Remember/reiterate & understand/summarize

| TASKS:<br>CONTEXT: | Work with/learn about Guidelines (GLs) | Engage with Stakeholder Analysis (SHA) | Learn and demonstrate Ethical Reasoning (KSAs 1-4) |
|---|---|---|---|
| Assumption (what if the assumption fails?) | -identify the GLs (from this list) that are relevant for ensuring | -if assumption/ approximation fail, or if application is not | -are the <given> GLs & SHA |
| Approximation (what if the approximation does not hold?) | assumption/ approximation holds or that application is justified/appropriate. -if SHA suggests harms to the public/public interest, which GLs are relevant? (Bloom's 1) What do the GLs suggest must be done? (Bloom's 2) | appropriate, do harms outweigh benefits (given a completed SHA) or not? -summarize the harms-benefits trade offs when assumption/ approximation fails or application is not justified | appropriate for <given vignette>? -is the ethical challenge correctly identified? -which <given> alternative is best supported? -given a case, fill in an ER table with KSA 1-4 as prompts |
| Application (is the application appropriate? what if it is not?) | | | |

*Learning outcomes midway through course (middle Blooms)*

Overall teaching objective: Engagement with the Ethical Reasoning process as dependent on stakeholders and how our decisions affect them, and how the target code/guidelines lead to the identification of ethical problems in the regular use of mathematics.

The *teaching* objectives are to reinforce the new schema for the student – that includes ethical reasoning and how their academic outputs represent their learning. Examples of learning outcomes fitting Stanford LO-writing criteria, include:

After this section of the course, students will:

1. Complete a case analysis utilizing all 6 Ethical Reasoning KSAs. (*learning how to complete a case analysis, introducing full range of* Ethical Reasoning *KSAs; structured communication about ethical mathematical practice*)
2. Evalute code/guidline elements for their relevance in specific instances where assumptions, approximations, or preconditions are not met or mathematical applications are not justified/appropriate across a variety of mathematical contexts. (*deepening familiarity with target code/guideline and how they function as a normative part of math)*
3. Identify harms and benefits to stakeholders whenever assumptions, approximations, or preconditions are not met or mathematical applications are not justified/appropriate across a variety of mathematical contexts. Complete the analysis of harms and benefits in a given context. (*deepening engagement with Stakeholder Analysis as a normative part of math; engaging in the analysis part of Stakeholder Analysis*)
4. Utilize stakeholder analysis and target code/guidelines in ethical reasoning about specific situations involving regular aspects of mathematics. (*practicing the* Ethical Reasoning *steps with a target code/guidelines and the Stakeholder*

*Analysis in typical mathematical applications; communicating about the ethical aspects of math*)

5. Describe how completing the tasks (below) represent evidence of their performance of each Ethical Reasoning KSA, including KSAs they have not learned/demonstrated yet. Compare earlier with later performance on each KSA. (*reinforcing engagement with all* Ethical Reasoning *KSAs and what evidence of their learning looks like*)

**Table 7.** Formula for creating an assignment (task) where learners must perform at Bloom's level 3-4: apply & analyze

| TASKS: CONTEXT: | work with Guidelines (GLs), construct and assess plans based on GLs | carry out and interpret Stakeholder Analysis (SHA) | learn and demonstrate full Ethical Reasoning KSA list (KSAs 1-6) |
|---|---|---|---|
| Assumption (what if the assumption fails?) | -determine which of a given set of GL elements) is the most strongly supportive of an ethical decision about assumption/ approximation failure or unjustified/ inappropriate application<br>- determine if an analysis of failure/unjustified application features the GLs correctly. | -complete your SHA for assumption/ approximation failure or inappropriate/ unjustified method<br>-outline harms-benefits tradeoff for your SHA | identify the GLs & complete the SHA appropriate for <given vignette>.<br>-is the ethical challenge correctly identified?<br>-which <given> alternative is best supported? |
| Approximation (what if the approximation does not hold?) | | | |
| Application (is the application appropriate? what if it is not?) | | | |

*Learning outcomes at end of course (high Blooms)*

Overall teaching objective: Demonstration of proficiency with the Ethical Reasoning process and communication about ethical aspects of mathematics, and their own learning.

The *teaching* objectives are to get learners to demonstrate how their new schema for mathematical practice includes ethical reasoning; and to reflect on how their academic outputs represent their learning. Examples of learning outcomes fitting Stanford LO-writing criteria, include:

After this section of the course, students will:

1. Complete a case analysis utilizing all six Ethical Reasoning KSAs, including a Stakeholder Analysis for each vignette. (*demonstrating proficient case analysis using full range of* Ethical Reasoning *KSAs; structured communication about ethical mathematical practice*)
2. Critically evaluate target code or guideline elements for their relevance in specific instances where assumptions, approximations, or preconditions are not met or mathematical applications are not justified/appropriate across a variety of mathematical contexts. (*deepening familiarity with target code/guidelines and how they function as a normative part of math)*
3. Analyze harms and benefits to stakeholders whenever assumptions, approximations, or preconditions are not met or mathematical applications are not justified/appropriate across a variety of mathematical contexts. Critically evaluate Stakeholder Analysis by others a given context, determining whether harms-benefits tradeoff has been appropriately articulated. (*deepening engagement with Stakeholder Analysis as a normative part of math; communicating about the impacts of decisionmaking on stakeholders in mathematical practice*)
4. Compare and contrast code or guidelines elements for their utility in decisionmaking in ethical reasoning about specific situations involving regular aspects of mathematics. (*critical thinking about codes/guidelines and what courses of action are best supported when ethical challenges arise in typical mathematical applications; communicating about the ethical aspects of math*)
5. Describe how completing the tasks (below) represent evidence of their performance of each Ethical Reasoning KSA. Compare earlier with later performance on each KSA. (*reinforcing engagement with all* Ethical Reasoning *KSAs and what evidence of their learning looks like*)

**Table 8.** Formula for creating an assignment (task) where learners must perform at Bloom's levels 5-6: create/apply/synthesize and evaluate/judge

| TASKS:<br>CONTEXT: | work with Guidelines (GLs), construct and assess plans based on GLs | carry out and interpret Stakeholder Analysis (SHA) | learn and demonstrate full Ethical Reasoning KSA list (KSAs 1-6) |
|---|---|---|---|
| Assumption (what if the assumption fails?) | Compare and contrast GLs that most and less strongly suggest a course of action given a vignette. | -complete SHA for vignettes<br>-evaluate a completed SHA for completeness<br>-critically assess a harms-benefits trade off (done by someone else) | -complete SHA and ER case analysis for vignette<br>- Explain the relationships between the MR-ER, and their work showing their level in the MR-ER, and |
| Approximation (what if the approximation does not hold?) | | | |
| Application (is the application | | | |
| appropriate? what if it is not?) | | | the curriculum in their degree/program. |

*Ethical Reasoning learning outcomes at end of curriculum/in Capstone (high Blooms)*

**Capstone course LOs presume that learners have been thinking about/working within the full Mastery Rubric for Ethical Reasoning: Students write 500 (max) word essays reflecting on the contents of these LOs.**
After the course, students will:
1. Describe the application of each Ethical Reasoning KSA –retroactively- to intellectual and academic experiences they have had (since high school and/or focused on current math course experiences);
2. Discuss the importance of the order of the Ethical Reasoning KSAs and in “real life”, using their own experiences or their case analyses from the term.
3. Rank order the Ethical Reasoning KSAs according to the student’s perception of their importance to competent and ethical practice (in their domain).
4. Explain the relationships between Ethical Reasoning and the standard “ethics” curriculum in their degree/program.
5. Describe 2-4 differences between a transcript record of their participation in this course and a "bingo card" representation of their performance/ capabilities/ achievements/plans on Mastery Rubric for Ethical Reasoning KSA (where they identify their performance levels, using evidence from their prior work, on each KSA).
6. Explain the evidence supporting their claim of their performance level on each KSA in the Mastery Rubric for Ethical Reasoning at the end of the course, and compare these with where they were at the start of the course.
7. Make and self-evaluate/revise a plan for moving to more advanced performance level descriptors s on two or more KSAs over the next (course, year, or other reasonable period).

As students craft these essays, they can utilize their own prior work -featuring their considerations of assumptions, approximations, and applications. They can also integrate other instances where failures in these three contexts did/did not result in benefits or harms to any particular stakeholders.

Students who are positioned to integrate the Ethical Reasoning KSAs into their thinking should be expected to be more advanced than those who were beginning to orient their thinking towards it, and our integrating course results demonstrated this. The 6$^{th}$ LO is specifically included as an exercise to evaluate (self-assess) the student’s evidence that they are not only “on the board”, but to begin to concretely describe themselves in terms of changes in their performance levels and capabilities on all relevant KSAs. This level of self-assessment and reflection will require relatively high Bloom’s level cognitive skills, and the teaching and practice in this integration course will require practice with feedback on self assessment (LO #7).

## 4. Discussion

While "ethics" is considered by some to be an important aspect of a mathematical or data science curriculum, once an individual enters professional practice - either in math/data science or in a job where they sometimes use math or data science - there is a much more universal expectation of behavior consistent with professional ethical practice: “Upon entry into practice, all professionals assume at least a tacit responsibility for the quality and integrity of their own work and that of colleagues. They also take on a responsibility to the larger public for the standards of practice associated with the profession.” (Golde & Walker, 2006: p. 10) This creates a challenge: how to integrate ethics content into quantitative courses without knowing the extent of responsibility to which professional codes of ethics any given student will ultimately face or take. Engagement with ethical reasoning can prepare students to assume this responsibility in any of its forms across their careers using mathematics, statistics, data science, and other quantitative methods and technologies. Moreover, leveraging a developmental trajectory over a course or curriculum can give instruction and practice with a wide variety of ethical reasoning situations and applications. Considering the integration of ethical reasoning beyond a single lesson or course can help to normalize ethical content, and reasoning, for quantitative practitioners.

APPENDIX

# ACM Code of Ethics and Professional Conduct

*Adopted by ACM Council 6/22/18.*

## Preamble

Computing professionals' actions change the world. To act responsibly, they should reflect upon the wider impacts of their work, consistently supporting the public good. The ACM Code of Ethics and Professional Conduct ("the Code") expresses the conscience of the profession.

The Code is designed to inspire and guide the ethical conduct of all computing professionals, including current and aspiring practitioners, instructors, students, influencers, and anyone who uses computing technology in an impactful way. Additionally, the Code serves as a basis for remediation when violations occur. The Code includes principles formulated as statements of responsibility, based on the understanding that the public good is always the primary consideration. Each principle is supplemented by guidelines, which provide explanations to assist computing professionals in understanding and applying the principle.

Section 1 outlines fundamental ethical principles that form the basis for the remainder of the Code. Section 2 addresses additional, more specific considerations of professional responsibility. Section 3 guides individuals who have a leadership role, whether in the workplace or in a volunteer professional capacity. Commitment to ethical conduct is required of every ACM member, and principles involving compliance with the Code are given in Section 4.

The Code as a whole is concerned with how fundamental ethical principles apply to a computing professional's conduct. The Code is not an algorithm for solving ethical problems; rather it serves as a basis for ethical decision-making. When thinking through a particular issue, a computing professional may find that multiple principles should be taken into account, and that different principles will have different relevance to the issue. Questions related to these kinds of issues can best be answered by thoughtful consideration of the fundamental ethical principles, understanding that the public good is the paramount consideration. The entire computing profession benefits when the ethical decision-making process is accountable to and transparent to all stakeholders. Open discussions about ethical issues promote this accountability and transparency.

## 1. GENERAL ETHICAL PRINCIPLES.

*A computing professional should…*

## 1.1 Contribute to society and to human well-being, acknowledging that all people are stakeholders in computing.

This principle, which concerns the quality of life of all people, affirms an obligation of computing professionals, both individually and collectively, to use their skills for the benefit of society, its members, and the environment surrounding them. This obligation includes promoting fundamental human rights and protecting each individual's right to autonomy. An essential aim of computing professionals is to minimize negative consequences of computing, including threats to health, safety, personal security, and privacy. When the interests of multiple groups conflict, the needs of those less advantaged should be given increased attention and priority.

Computing professionals should consider whether the results of their efforts will respect diversity, will be used in socially responsible ways, will meet social needs, and will be broadly accessible. They are encouraged to actively contribute to society by engaging in pro bono or volunteer work that benefits the public good.

In addition to a safe social environment, human well-being requires a safe natural environment. Therefore, computing professionals should promote environmental sustainability both locally and globally.

## 1.2 Avoid harm.

In this document, “harm” means negative consequences, especially when those consequences are significant and unjust. Examples of harm include unjustified physical or mental injury, unjustified destruction or disclosure of information, and unjustified damage to property, reputation, and the environment. This list is not exhaustive.

Well-intended actions, including those that accomplish assigned duties, may lead to harm. When that harm is unintended, those responsible are obliged to undo or mitigate the harm as much as possible. Avoiding harm begins with careful consideration of potential impacts on all those affected by decisions. When harm is an intentional part of the system, those responsible are obligated to ensure that the harm is ethically justified. In either case, ensure that all harm is minimized.

To minimize the possibility of indirectly or unintentionally harming others, computing professionals should follow generally accepted best practices unless there is a compelling ethical reason to do otherwise. Additionally, the consequences of data aggregation and emergent properties of systems should be carefully analyzed. Those involved with pervasive or infrastructure systems should also consider Principle 3.7.

A computing professional has an additional obligation to report any signs of system risks that might result in harm. If leaders do not act to curtail or mitigate such risks, it may be necessary to “blow the whistle” to reduce potential harm. However, capricious or misguided reporting of risks can itself be harmful. Before reporting risks, a computing professional should carefully assess relevant aspects of the situation.

## 1.3 Be honest and trustworthy.

Honesty is an essential component of trustworthiness. A computing professional should be transparent and provide full disclosure of all pertinent system capabilities, limitations, and potential problems to the appropriate parties. Making deliberately false or misleading claims, fabricating or falsifying data, offering or accepting bribes, and other dishonest conduct are violations of the Code.

Computing professionals should be honest about their qualifications, and about any limitations in their competence to complete a task. Computing professionals should be forthright about any circumstances that might lead to either real or perceived conflicts of interest or otherwise tend to undermine the independence of their judgment. Furthermore, commitments should be honored.

Computing professionals should not misrepresent an organization’s policies or procedures, and should not speak on behalf of an organization unless authorized to do so.

## 1.4 Be fair and take action not to discriminate.

The values of equality, tolerance, respect for others, and justice govern this principle. Fairness requires that even careful decision processes provide some avenue for redress of grievances.

Computing professionals should foster fair participation of all people, including those of underrepresented groups. Prejudicial discrimination on the basis of age, color, disability, ethnicity, family status, gender identity, labor union membership, military status, nationality, race, religion or belief, sex, sexual orientation, or any other inappropriate factor is an explicit violation of the Code. Harassment, including sexual harassment, bullying, and other abuses of power and authority, is a form of discrimination that, amongst other harms, limits fair access to the virtual and physical spaces where such harassment takes place.

The use of information and technology may cause new, or enhance existing, inequities. Technologies and practices should be as inclusive and accessible as possible and computing professionals should take action to avoid creating systems or technologies that disenfranchise or oppress people. Failure to design for inclusiveness and accessibility may constitute unfair discrimination.

## 1.5 Respect the work required to produce new ideas, inventions, creative works, and computing artifacts.

Developing new ideas, inventions, creative works, and computing artifacts creates value for society, and those who expend this effort should expect to gain value from their work. Computing professionals should therefore credit the creators of ideas, inventions, work, and artifacts, and respect copyrights, patents, trade secrets, license agreements, and other methods of protecting authors' works.

Both custom and the law recognize that some exceptions to a creator's control of a work are necessary for the public good. Computing professionals should not unduly oppose reasonable uses of their intellectual works. Efforts to help others by contributing time and energy to projects that help society illustrate a positive aspect of this principle. Such efforts include free and open source software and work put into the public domain. Computing professionals should not claim private ownership of work that they or others have shared as public resources.

## 1.6 Respect privacy.

The responsibility of respecting privacy applies to computing professionals in a particularly profound way. Technology enables the collection, monitoring, and exchange of personal information quickly, inexpensively, and often without the knowledge of the people affected. Therefore, a computing professional should become conversant in the various definitions and forms of privacy and should understand the rights and responsibilities associated with the collection and use of personal information.

Computing professionals should only use personal information for legitimate ends and without violating the rights of individuals and groups. This requires taking precautions to prevent re- identification of anonymized data or unauthorized data collection, ensuring the accuracy of data, understanding the provenance of the data, and protecting it from unauthorized access and accidental disclosure. Computing professionals should establish transparent policies and procedures that allow individuals to understand what data is being collected and how it is being used, to give informed consent for automatic data collection, and to review, obtain, correct inaccuracies in, and delete their personal data.

Only the minimum amount of personal information necessary should be collected in a system. The retention and disposal periods for that information should be clearly defined, enforced, and communicated to data subjects. Personal information gathered for a specific purpose should not be used for other purposes without the person's consent. Merged data collections can compromise privacy features present in the original collections. Therefore, computing professionals should take special care for privacy when merging data collections.

### 1.7 Honor confidentiality.

Computing professionals are often entrusted with confidential information such as trade secrets, client data, nonpublic business strategies, financial information, research data, pre-publication scholarly articles, and patent applications. Computing professionals should protect confidentiality except in cases where it is evidence of the violation of law, of organizational regulations, or of the Code. In these cases, the nature or contents of that information should not be disclosed except to appropriate authorities. A computing professional should consider thoughtfully whether such disclosures are consistent with the Code.

## 2. PROFESSIONAL RESPONSIBILITIES.

*A computing professional should…*

### 2.1 Strive to achieve high quality in both the processes and products of professional work.

Computing professionals should insist on and support high quality work from themselves and from colleagues. The dignity of employers, employees, colleagues, clients, users, and anyone else affected either directly or indirectly by the work should be respected throughout the process. Computing professionals should respect the right of those involved to transparent communication about the project. Professionals should be cognizant of any serious negative consequences affecting any stakeholder that may result from poor quality work and should resist inducements to neglect this responsibility.

### 2.2 Maintain high standards of professional competence, conduct, and ethical practice.

High quality computing depends on individuals and teams who take personal and group responsibility for acquiring and maintaining professional competence. Professional competence starts with technical knowledge and with awareness of the social context in which their work may be deployed. Professional competence also requires skill in communication, in reflective analysis, and in recognizing and navigating ethical challenges. Upgrading skills should be an ongoing process and might include independent study, attending conferences or seminars, and other informal or formal education. Professional organizations and employers should encourage and facilitate these activities.

## 2.3 Know and respect existing rules pertaining to professional work.

“Rules” here include local, regional, national, and international laws and regulations, as well as any policies and procedures of the organizations to which the professional belongs. Computing professionals must abide by these rules unless there is a compelling ethical justification to do otherwise. Rules that are judged unethical should be challenged. A rule may be unethical when it has an inadequate moral basis or causes recognizable harm. A computing professional should consider challenging the rule through existing channels before violating the rule. A computing professional who decides to violate a rule because it is unethical, or for any other reason, must consider potential consequences and accept responsibility for that action.

## 2.4 Accept and provide appropriate professional review.

High quality professional work in computing depends on professional review at all stages. Whenever appropriate, computing professionals should seek and utilize peer and stakeholder review. Computing professionals should also provide constructive, critical reviews of others’ work.

## 2.5 Give comprehensive and thorough evaluations of computer systems and their impacts, including analysis of possible risks.

Computing professionals are in a position of trust, and therefore have a special responsibility to provide objective, credible evaluations and testimony to employers, employees, clients, users, and the public. Computing professionals should strive to be perceptive, thorough, and objective when evaluating, recommending, and presenting system descriptions and alternatives. Extraordinary care should be taken to identify and mitigate potential risks in machine learning systems. A system for which future risks cannot be reliably predicted requires frequent reassessment of risk as the system evolves in use, or it should not be deployed. Any issues that might result in major risk must be reported to appropriate parties.

## 2.6 Perform work only in areas of competence.

A computing professional is responsible for evaluating potential work assignments. This includes evaluating the work’s feasibility and advisability, and making a judgment about whether the work assignment is within the professional’s areas of competence. If at any time before or during the work assignment the professional identifies a lack of a necessary expertise, they must disclose this to the employer or client. The client or employer may decide to pursue the assignment with the professional after additional time to acquire the necessary competencies, to pursue the assignment with someone else who has the required expertise, or to forgo the assignment. A computing professional's ethical judgment should be the final guide in deciding whether to work on the assignment.

### 2.7 Foster public awareness and understanding of computing, related technologies, and their consequences.

As appropriate to the context and one's abilities, computing professionals should share technical knowledge with the public, foster awareness of computing, and encourage understanding of computing. These communications with the public should be clear, respectful, and welcoming. Important issues include the impacts of computer systems, their limitations, their vulnerabilities, and the opportunities that they present. Additionally, a computing professional should respectfully address inaccurate or misleading information related to computing.

### 2.8 Access computing and communication resources only when authorized or when compelled by the public good.

Individuals and organizations have the right to restrict access to their systems and data so long as the restrictions are consistent with other principles in the Code. Consequently, computing professionals should not access another's computer system, software, or data without a reasonable belief that such an action would be authorized or a compelling belief that it is consistent with the public good. A system being publicly accessible is not sufficient grounds on its own to imply authorization. Under exceptional circumstances a computing professional may use unauthorized access to disrupt or inhibit the functioning of malicious systems; extraordinary precautions must be taken in these instances to avoid harm to others.

### 2.9 Design and implement systems that are robustly and usably secure.

Breaches of computer security cause harm. Robust security should be a primary consideration when designing and implementing systems. Computing professionals should perform due diligence to ensure the system functions as intended, and take appropriate action to secure resources against accidental and intentional misuse, modification, and denial of service. As threats can arise and change after a system is deployed, computing professionals should integrate mitigation techniques and policies, such as monitoring, patching, and vulnerability reporting. Computing professionals should also take steps to ensure parties affected by data breaches are notified in a timely and clear manner, providing appropriate guidance and remediation.

To ensure the system achieves its intended purpose, security features should be designed to be as intuitive and easy to use as possible. Computing professionals should discourage security precautions that are too confusing, are situationally inappropriate, or otherwise inhibit legitimate use.

In cases where misuse or harm are predictable or unavoidable, the best option may be to not implement the system.

# 3. PROFESSIONAL LEADERSHIP PRINCIPLES.

Leadership may either be a formal designation or arise informally from influence over others. In this section, "leader" means any member of an organization or group who has influence, educational responsibilities, or managerial responsibilities. While these principles apply to all computing professionals, leaders bear a heightened responsibility to uphold and promote them, both within and through their organizations.

*A computing professional, especially one acting as a leader, should…*

## 3.1 Ensure that the public good is the central concern during all professional computing work.

People—including users, customers, colleagues, and others affected directly or indirectly— should always be the central concern in computing. The public good should always be an explicit consideration when evaluating tasks associated with research, requirements analysis, design, implementation, testing, validation, deployment, maintenance, retirement, and disposal. Computing professionals should keep this focus no matter which methodologies or techniques they use in their practice.

## 3.2 Articulate, encourage acceptance of, and evaluate fulfillment of social responsibilities by members of the organization or group.

Technical organizations and groups affect broader society, and their leaders should accept the associated responsibilities. Organizations—through procedures and attitudes oriented toward quality, transparency, and the welfare of society—reduce harm to the public and raise awareness of the influence of technology in our lives. Therefore, leaders should encourage full participation of computing professionals in meeting relevant social responsibilities and discourage tendencies to do otherwise.

### 3.3 Manage personnel and resources to enhance the quality of working life.

Leaders should ensure that they enhance, not degrade, the quality of working life. Leaders should consider the personal and professional development, accessibility requirements, physical safety, psychological well-being, and human dignity of all workers. Appropriate human-computer ergonomic standards should be used in the workplace.

### 3.4 Articulate, apply, and support policies and processes that reflect the principles of the Code.

Leaders should pursue clearly defined organizational policies that are consistent with the Code and effectively communicate them to relevant stakeholders. In addition, leaders should encourage and reward compliance with those policies, and take appropriate action when policies are violated. Designing or implementing processes that deliberately or negligently violate, or tend to enable the violation of, the Code's principles is ethically unacceptable.

### 3.5 Create opportunities for members of the organization or group to grow as professionals.

Educational opportunities are essential for all organization and group members. Leaders should ensure that opportunities are available to computing professionals to help them improve their knowledge and skills in professionalism, in the practice of ethics, and in their technical specialties. These opportunities should include experiences that familiarize computing professionals with the consequences and limitations of particular types of systems. Computing professionals should be fully aware of the dangers of oversimplified approaches, the improbability of anticipating every possible operating condition, the inevitability of software errors, the interactions of systems and their contexts, and other issues related to the complexity of their profession—and thus be confident in taking on responsibilities for the work that they do.

### 3.6 Use care when modifying or retiring systems.

Interface changes, the removal of features, and even software updates have an impact on the productivity of users and the quality of their work. Leaders should take care when changing or discontinuing support for system features on which people still depend. Leaders should thoroughly investigate viable alternatives to removing support for a legacy system. If these alternatives are unacceptably risky or impractical, the developer should assist stakeholders' graceful migration from the system to an alternative. Users should be notified of the risks of continued use of the unsupported system long before support ends. Computing professionals should assist system users in monitoring the operational viability of their computing systems, and help them understand that timely replacement of inappropriate or outdated features or entire systems may be needed.

### 3.7 Recognize and take special care of systems that become integrated into the infrastructure of society.

Even the simplest computer systems have the potential to impact all aspects of society when integrated with everyday activities such as commerce, travel, government, healthcare, and education. When organizations and groups develop systems that become an important part of the infrastructure of society, their leaders have an added responsibility to be good stewards of these systems. Part of that stewardship requires establishing policies for fair system access, including for those who may have been excluded. That stewardship also requires that computing professionals monitor the level of integration of their systems into the infrastructure of society. As the level of adoption changes, the ethical responsibilities of the organization or group are likely to change as well. Continual monitoring of how society is using a system will allow the organization or group to remain consistent with their ethical obligations outlined in the Code. When appropriate standards of care do not exist, computing professionals have a duty to ensure they are developed.

## 4. COMPLIANCE WITH THE CODE.

*A computing professional should…*

### 4.1 Uphold, promote, and respect the principles of the Code.

The future of computing depends on both technical and ethical excellence. Computing professionals should adhere to the principles of the Code and contribute to improving them. Computing professionals who recognize breaches of the Code should take actions to resolve the ethical issues they recognize, including, when reasonable, expressing their concern to the person or persons thought to be violating the Code.

### 4.2 Treat violations of the Code as inconsistent with membership in the ACM.

Each ACM member should encourage and support adherence by all computing professionals regardless of ACM membership. ACM members who recognize a breach of the Code should consider reporting the violation to the ACM, which may result in remedial action as specified in the ACM's Code of Ethics and Professional Conduct Enforcement Policy.

---

*The Code and guidelines were developed by the ACM Code 2018 Task Force: Executive Committee Don Gotterbarn (Chair), Bo Brinkman, Catherine Flick, Michael S Kirkpatrick, Keith Miller, Kate Varansky, and Marty J Wolf. Members: Eve Anderson, Ron Anderson, Amy Bruckman, Karla Carter, Michael Davis, Penny Duquenoy, Jeremy Epstein, Kai Kimppa, Lorraine Kisselburgh, Shrawan Kumar, Andrew McGettrick, Natasa Milic-Frayling, Denise Oram, Simon Rogerson, David Shama, Janice Sipior, Eugene Spafford, and Les Waguespack. The Task Force was organized by the ACM Committee on Professional Ethics. Significant contributions to the Code were also made by the broader international ACM membership. This Code and its guidelines were adopted by the ACM Council on June 22nd, 2018.*

## Ethical Guidelines for Statistical Practice

***Prepared by the Committee on Professional Ethics of the American Statistical Association***

**February 2022**

### PURPOSE OF THE GUIDELINES:

The American Statistical Association's Ethical Guidelines for Statistical Practice are intended to help statistical practitioners make decisions ethically. In these Guidelines, “statistical practice” includes activities such as: designing the collection of, summarizing, processing, analyzing, interpreting, or presenting, data; as well as model or algorithm development and deployment. Throughout these Guidelines, the term "statistical practitioner" includes all those who engage in statistical practice, regardless of job title, profession, level, or field of degree. The Guidelines are intended for individuals, but these principles are also relevant to organizations that engage in statistical practice.

The Ethical Guidelines aim to promote accountability by informing those who rely on any aspects of statistical practice of the standards that they should expect. Society benefits from informed judgments supported by ethical statistical practice. All statistical practitioners are expected to follow these Guidelines and to encourage others to do the same.

In some situations, Guideline principles may require balancing of competing interests. If an unexpected ethical challenge arises, the ethical practitioner seeks guidance, not exceptions, in the Guidelines. To justify unethical behaviors, or to exploit gaps in the Guidelines, is unprofessional, and inconsistent with these Guidelines.

### PRINCIPLE A: Professional Integrity and Accountability

Professional integrity and accountability require taking responsibility for one's work. Ethical statistical practice supports valid and prudent decision making with appropriate methodology. The ethical statistical practitioner represents their capabilities and activities honestly, and treats others with respect.

**The ethical statistical practitioner:**

1. Takes responsibility for evaluating potential tasks, assessing whether they have (or can attain) sufficient competence to execute each task, and that the work and timeline are feasible. Does not solicit or deliver work for which they are not qualified, or that they would not be willing to have peer reviewed.
2. Uses methodology and data that are valid, relevant, and appropriate, without favoritism or prejudice, and in a manner intended to produce valid, interpretable, and reproducible results.
3. Does not knowingly conduct statistical practices that exploit vulnerable populations or create or perpetuate unfair outcomes.
4. Opposes efforts to predetermine or influence the results of statistical practices, and resists pressure to selectively interpret data.
5. Accepts full responsibility for their own work; does not take credit for the work of others; and gives credit to those who contribute. Respects and acknowledges the intellectual property of others.
6. Strives to follow, and encourages all collaborators to follow, an established protocol for authorship. Advocates for recognition commensurate with each person's contribution to the work. Recognizes that inclusion as an author does imply, while acknowledgement may imply, endorsement of the work.
7. Discloses conflicts of interest, financial and otherwise, and manages or resolves them according to established policies, regulations, and laws.
8. Promotes the dignity and fair treatment of all people. Neither engages in nor condones discrimination based on personal characteristics. Respects personal boundaries in interactions and avoids harassment including sexual harassment, bullying, and other abuses of power or authority.
9. Takes appropriate action when aware of deviations from these Guidelines by others.
10. Acquires and maintains competence through upgrading of skills as needed to maintain a high standard of practice.
11. Follows applicable policies, regulations, and laws relating to their professional work, unless there is a compelling ethical justification to do otherwise.
12. Upholds, respects, and promotes these Guidelines. Those who teach, train, or mentor in statistical practice have a special obligation to promote behavior that is consistent with these Guidelines.

**PRINCIPLE B: Integrity of Data and Methods**

The ethical statistical practitioner seeks to understand and mitigate known or suspected limitations, defects, or biases in the data or methods and communicates potential impacts on the interpretation, conclusions, recommendations, decisions, or other results of statistical practices.

**The ethical statistical practitioner:**

1. Communicates data sources and fitness for use, including data generation and collection processes and known biases. Discloses and manages any conflicts of interest relating to the data sources. Communicates data processing and transformation procedures, including missing data handling.
2. Is transparent about assumptions made in the execution and interpretation of statistical practices including methods used, limitations, possible sources of error, and algorithmic biases. Conveys results or applications of statistical practices in ways that are honest and meaningful.
3. Communicates the stated purpose and the intended use of statistical practices. Is transparent regarding a priori versus post hoc objectives and planned versus unplanned statistical practices. Discloses when multiple comparisons are conducted, and any relevant adjustments.
4. Meets obligations to share the data used in the statistical practices, for example, for peer review and replication, as allowable. Respects expectations of data contributors when using or sharing data. Exercises due caution to protect proprietary and confidential data, including all data that might inappropriately harm data subjects.
5. Strives to promptly correct substantive errors discovered after publication or implementation. As appropriate, disseminates the correction publicly and/or to others relying on the results.
6. For models and algorithms designed to inform or implement decisions repeatedly, develops and/or implements plans to validate assumptions and assess performance over time, as needed. Considers criteria and mitigation plans for model or algorithm failure and retirement.
7. Explores and describes the effect of variation in human characteristics and groups on statistical practice when feasible and relevant.

## PRINCIPLE C: Responsibilities to Stakeholders

Those who fund, contribute to, use, or are affected by statistical practices are considered stakeholders. The ethical statistical practitioner respects the interests of stakeholders while practicing in compliance with these Guidelines.

**The ethical statistical practitioner:**

1. Seeks to establish what stakeholders hope to obtain from any specific project. Strives to obtain sufficient subject-matter knowledge to conduct meaningful and relevant statistical practice.
2. Regardless of personal or institutional interests or external pressures, does not use statistical practices to mislead any stakeholder.
3. Uses practices appropriate to exploratory and confirmatory phases of a project, differentiating findings from each so the stakeholders can understand and apply the results.
4. Informs stakeholders of the potential limitations on use and re-use of statistical practices in different contexts and offers guidance and alternatives, where appropriate, about scope, cost, and precision considerations that affect the utility of the statistical practice.
5. Explains any expected adverse consequences from failing to follow through on an agreed-upon sampling or analytic plan.
6. Strives to make new methodological knowledge widely available to provide benefits to society at large. Presents relevant findings, when possible, to advance public knowledge.

7. Understands and conforms to confidentiality requirements for data collection, release, and dissemination and any restrictions on its use established by the data provider (to the extent legally required). Protects the use and disclosure of data accordingly. Safeguards privileged information of the employer, client, or funder.
8. Prioritizes both scientific integrity and the principles outlined in these Guidelines when interests are in conflict.

## PRINCIPLE D: Responsibilities to Research Subjects, Data Subjects, or those directly affected by statistical practices

The ethical statistical practitioner does not misuse or condone the misuse of data. They protect and respect the rights and interests of human and animal subjects. These responsibilities extend to those who will be directly affected by statistical practices.

**The ethical statistical practitioner:**

1. Keeps informed about and adheres to applicable rules, approvals, and guidelines for the protection and welfare of human and animal subjects. Knows when work requires ethical review and oversight.[3]
2. Makes informed recommendations for sample size and statistical practice methodology in order to avoid the use of excessive or inadequate numbers of subjects and excessive risk to subjects
3. For animal studies, seeks to leverage statistical practice to reduce the number of animals used, refine experiments to increase the humane treatment of animals, and replace animal use where possible.
4. Protects people's privacy and the confidentiality of data concerning them, whether obtained from the individuals directly, other persons, or existing records. Knows and adheres to applicable rules, consents, and guidelines to protect private information.
5. Uses data only as permitted by data subjects' consent when applicable or considering their interests and welfare when consent is not required. This includes primary and secondary uses, use of repurposed data, sharing data, and linking data with additional data sets.
6. Considers the impact of statistical practice on society, groups, and individuals. Recognizes that statistical practice could adversely affect groups or the public perception of groups, including marginalized groups. Considers approaches to minimize negative impacts in applications or in framing results in reporting.
7. Refrains from collecting or using more data than is necessary. Uses confidential information only when permitted and only to the extent necessary. Seeks to minimize the risk of re-identification when sharing de-identified data or results where there is an expectation of confidentiality. Explains any impact of de-identification on accuracy of results.
8. To maximize contributions of data subjects, considers how best to use available data sources for exploration, training, testing, validation, or replication as needed for the application. The ethical statistical practitioner appropriately discloses how the data is used for these purposes and any limitations.

---

[3] Examples of ethical review and oversight include an Institutional Review Board (IRB), an Institutional Animal Care and Use Committee (IACUC), or a compliance assessment.

9. Knows the legal limitations on privacy and confidentiality assurances and does not over-promise or assume legal privacy and confidentiality protections where they may not apply.
10. Understands the provenance of the data, including origins, revisions, and any restrictions on usage, and fitness for use prior to conducting statistical practices.
11. Does not conduct statistical practice that could reasonably be interpreted by subjects as sanctioning a violation of their rights. Seeks to use statistical practices to promote the just and impartial treatment of all individuals.

**PRINCIPLE E: Responsibilities to members of multidisciplinary teams**

Statistical practice is often conducted in teams made up of professionals with different professional standards. The statistical practitioner must know how to work ethically in this environment.

**The ethical statistical practitioner:**

1. Recognizes and respects that other professions may have different ethical standards and obligations. Dissonance in ethics may still arise even if all members feel that they are working towards the same goal. It is essential to have a respectful exchange of views.
2. Prioritizes these Guidelines for the conduct of statistical practice in cases where ethical guidelines conflict.
3. Ensures that all communications regarding statistical practices are consistent with these Guidelines. Promotes transparency in all statistical practices.
4. Avoids compromising validity for expediency. Regardless of pressure on or within the team, does not use inappropriate statistical practices.

**PRINCIPLE F: Responsibilities to Fellow Statistical Practitioners and the Profession**

Statistical practices occur in a wide range of contexts. Irrespective of job title and training, those who practice statistics have a responsibility to treat statistical practitioners, and the profession, with respect. Responsibilities to other practitioners and the profession include honest communication and engagement that can strengthen the work of others and the profession.

**The ethical statistical practitioner:**

1. Recognizes that statistical practitioners may have different expertise and experiences, which may lead to divergent judgments about statistical practices and results. Constructive discourse with mutual respect focuses on scientific principles and methodology and not personal attributes.
2. Helps strengthen, and does not undermine, the work of others through appropriate peer review or consultation. Provides feedback or advice that is impartial, constructive, and objective.
3. Takes full responsibility for their contributions as instructors, mentors, and supervisors of statistical practice by ensuring their best teaching and advising -- regardless of an academic or non-academic setting -- to ensure that developing practitioners are guided effectively as they learn and grow in their careers.

4. Promotes reproducibility and replication, whether results are “significant” or not, by sharing data, methods, and documentation to the extent possible.
5. Serves as an ambassador for statistical practice by promoting thoughtful choices about data acquisition, analytic procedures, and data structures among non-practitioners and students. Instills appreciation for the concepts and methods of statistical practice.

## PRINCIPLE G: Responsibilities of Leaders, Supervisors, and Mentors in Statistical Practice

Statistical practitioners leading, supervising, and/or mentoring people in statistical practice have specific obligations to follow and promote these Ethical Guidelines. Their support for – and insistence on – ethical statistical practice are essential for the integrity of the practice and profession of statistics as well as the practitioners themselves.

**Those leading, supervising, or mentoring statistical practitioners are expected to**:

1. Ensure appropriate statistical practice that is consistent with these Guidelines. Protect the statistical practitioners who comply with these Guidelines, and advocate for a culture that supports ethical statistical practice.
2. Promote a respectful, safe, and productive work environment. Encourage constructive engagement to improve statistical practice.
3. Identify and/or create opportunities for team members/mentees to develop professionally and maintain their proficiency.
4. Advocate for appropriate, timely, inclusion and participation of statistical practitioners as contributors/collaborators. Promote appropriate recognition of the contributions of statistical practitioners, including authorship if applicable.
5. Establish a culture that values validation of assumptions, and assessment of model/algorithm performance over time and across relevant subgroups, as needed. Communicate with relevant stakeholders regarding model or algorithm maintenance, failure, or actual or proposed modifications.

## PRINCIPLE H: Responsibilities Regarding Potential Misconduct

The ethical statistical practitioner understands that questions may arise concerning potential misconduct related to statistical, scientific, or professional practice. At times, a practitioner may accuse someone of misconduct, or be accused by others. At other times, a practitioner may be involved in the investigation of others’ behavior. Allegations of misconduct may arise within different institutions with different standards and potentially different outcomes. The elements that follow relate specifically to allegations of statistical, scientific, and professional misconduct.

**The ethical statistical practitioner:**

1. Knows the definitions of, and procedures relating to, misconduct in their institutional setting. Seeks to clarify facts and intent before alleging misconduct by others. Recognizes that differences of opinion and honest error do not constitute unethical behavior.

2. Avoids condoning or appearing to condone statistical, scientific, or professional misconduct. Encourages other practitioners to avoid misconduct or the appearance of misconduct.
3. Does not make allegations that are poorly founded, or intended to intimidate. Recognizes such allegations as potential ethics violations.
4. Lodges complaints of misconduct discreetly and to the relevant institutional body. Does not act on allegations of misconduct without appropriate institutional referral, including those allegations originating from social media accounts or email listservs.
5. Insists upon a transparent and fair process to adjudicate claims of misconduct. Maintains confidentiality when participating in an investigation. Discloses the investigation results honestly to appropriate parties and stakeholders once they are available.
6. Refuses to publicly question or discredit the reputation of a person based on a specific accusation of misconduct while due process continues to unfold.
7. Following an investigation of misconduct, supports the efforts of all parties involved to resume their careers in as normal a manner as possible, consistent with the outcome of the investigation.
8. Avoids, and acts to discourage, retaliation against or damage to the employability of those who responsibly call attention to possible misconduct.

## APPENDIX
## Responsibilities of organizations/institutions

Whenever organizations and institutions design the collection of, summarize, process, analyze, interpret, or present, data; or develop and/or deploy models or algorithms, they have responsibilities to use statistical practice in ways that are consistent with these Guidelines, as well as promote ethical statistical practice.

**Organizations and institutions engage in, and promote, ethical statistical practice by**:

1. Expecting and encouraging all employees and vendors who conduct statistical practice to adhere to these Guidelines. Promoting a workplace where the ethical practitioner may apply the Guidelines without being intimidated or coerced. Protecting statistical practitioners who comply with these Guidelines.
2. Engaging competent personnel to conduct statistical practice, and promote a productive work environment.
3. Promoting the professional development and maintenance of proficiency for employed statistical practitioners.
4. Supporting statistical practice that is objective and transparent. Not allowing organizational objectives or expectations to encourage unethical statistical practice by its employees.
5. Recognizing that the inclusion of statistical practitioners as authors, or acknowledgement of their contributions to projects or publications, requires their explicit permission because it may imply endorsement of the work.
6. Avoiding statistical practices that exploit vulnerable populations or create or perpetuate discrimination or unjust outcomes. Considering both scientific validity and impact on societal and human well-being that results from the organization’s statistical practice.
7. Using professional qualifications and contributions as the basis for decisions regarding statistical practitioners' hiring, firing, promotion, work assignments, publications and presentations, candidacy for offices and awards, funding or approval of research, and other professional matters.

**Those in leadership, supervisory, or managerial positions who oversee statistical practitioners promote ethical statistical practice by following Principle G and:**

8. Recognizing that it is contrary to these Guidelines to report or follow only those results that conform to expectations without explicitly acknowledging competing findings and the basis for choices regarding which results to report, use, and/or cite.
9. Recognizing that the results of valid statistical studies cannot be guaranteed to conform to the expectations or desires of those commissioning the study or employing/supervising the statistical practitioner(s).
10. Objectively, accurately, and efficiently communicating a team's or practitioners' statistical work throughout the organization.
11. In cases where ethical issues are raised, representing them fairly within the organization's leadership team.
12. Managing resources and organizational strategy to direct teams of statistical practitioners along the most productive lines in light of the ethical standards contained in these Guidelines.

Ethical Proto-Guidelines for Mathematical Practice (Tractenberg et al. 2023; used with permission)

| **The ethical mathematics practitioner...** |
|---|
| **IN GENERAL** |
| 1. Is honest about their qualification to complete work they accept; articulates any limitation of expertise, and consults others when necessary or in doubt. They accept full responsibility for their professional performance and practice. |
| 2. Treats others with respect. Promotes the equal dignity and fair treatment of all people, and neither engages in nor condones discrimination based on personal characteristics. Respects personal boundaries in interactions, and avoids harassment, including sexual harassment; bullying; and other abuses of power or authority. Takes appropriate action when aware of disrespectful behaviors by others. |
| 3. Accepts full responsibility for their own work; does not take credit for the work of others; and gives credit to those who contribute. Respects and acknowledges the intellectual property of others. |
| 4. Should be forthright about any circumstances that might lead to either real or perceived conflicts of interest or otherwise tend to undermine the independence of their judgment. Discloses conflicts of interest, financial and otherwise, and manages or resolves them according to established (institutional/regional/local) rules and laws. |
| 5. Recognizes any mathematical descriptions of groups may carry risks of stereotypes and stigmatization. Practitioners should contemplate, and be sensitive to, the manner in which information in their work across education, research, |
| public policy, and in the public in general, is framed to avoid disproportionate harm to vulnerable groups. |
| 6. Avoids condoning or appearing to condone mathematical, scientific, or professional misconduct. Takes appropriate action when aware of unethical conduct by others. |
| 7. Avoids, and acts to discourage, retaliation against or damage to the employability of those who responsibly call attention to possible mathematical error or to scientific or other misconduct. |
| 8. Is informed about applicable laws, policies, rules, and guidelines; follows these unless there is a compelling ethical reason to do otherwise. |
| 9. Must know how to work ethically in collaborative environment. When conducting their work in conjunction with other professions, must continue to abide by mathematicians' responsibilities, as well as any guidelines of the other professions. When there is a conflict or an absence in the partner profession's guidelines, the mathematical practitioners' responsibilities should be followed. |
| 10. Respects others, and promotes justice and inclusiveness, in all work. Fosters fair participation of all people. Avoids and mitigates bias and prejudice. Does nothing to limit fair access. |
| 11. Opposes marginalization of people on the basis of human differences. Strives to resist institutional confirmation bias and systematic injustice. |
| 12. Minimizes the possibility of harming others; whether directly or indirectly, intentionally or unintentionally. |
| **AS A MEMBER OF THE PROFESSION** |
| 13. Strives to make new mathematical knowledge as widely available as is feasible. |
| 14. Maintains high standards of professional competence, conduct, and ethical practices. |
| 15. Recognizes that if they engage in mathematics practice, they do so in a social and cultural context, acknowledging that all people are stakeholders in mathematics. |
| 16. In reviews, considers the potential for unjust or inequitable implications of the proposal or work. |
| 17. Understands the differences between questionable mathematical, scientific, or professional practices and practices that constitute misconduct. The ethical mathematics practitioner avoids all of the above and knows how each should be handled. |
| 18. Helps strengthen the work of others through appropriate peer review; assesses methods, not individuals. Strives to complete review assignments thoroughly, thoughtfully, and promptly. |
| 19. Avoids and addresses exclusionary practices in hiring, teaching, and recruitinWhen assessing or evaluating mathematics practitioners or their work, uses relevant subject matter-specific qualifications. Uses qualifications, performance, and contributions as the basis for decisions regarding mathematical practitioners of all levels. |
| 20. Upholds, promotes, and respects the ethical responsibilities of the mathematics community. |
| 21. Accepts their accountability to build an inclusive mathematics community that values its members. |
| 22. When involved in advising graduate students, should fully inform them about the employment prospects they may face upon completion of their degrees. |
| **IN THEIR SCHOLARSHIP** |
| 23. Strives to support and achieve quality work in both the process and products of professional work. Works in a manner intended to produce valid, interpretable, and when applicable, reproducible results. |
| 24. Identifies and mitigates any efforts to predetermine or influence the results or outcomes of mathematical practices; resists pressure to solve unethical problems/support predetermined outcomes. |
| 25. Strives to follow, and encourages all collaborators to follow, an established protocol for authorship. |
| 26. Is candid about any known or suspected limitations, assumptions, or biases when working with methods, models, or data. Objective and valid interpretation of the results requires that the underlying analysis recognizes and acknowledges the degree of reliability and integrity of the method, model, or data. |
| 27. Assesses, and is transparent about, the origin and source of the tools and methods they use, including prior results and data. Practitioners, when possible, acknowledge and disclose the origin of the problems they are solving and the interests that their work is intended to serve. |
| 28. Strives to promptly correct any errors discovered while producing the final report or after publication. As appropriate, disseminates the correction publicly or to others relying on the results. |
| 29. Understands and conforms to confidentiality requirements of data collection, release, and dissemination and any restrictions on its use established by the data provider (to the extent legally required), protecting use and disclosure of data accordingly. |
| 30. Strives to ensure that data sources, choice of methods, and applications do not create or perpetuate social biases or discrimination. Seeks to avoid confirmation bias. |
| 31. Avoids plagiarism. The knowing presentation of another person's mathematical discovery as one's own constitutes plagiarism and is a serious violation of professional ethics. Plagiarism may occur for any type of work, whether written or oral and whether published or not. |
| 32. Promotes sharing of data, methods, scholarship as much as possible and as appropriate without compromising propriety. |
| 33. Recognizes the inclusion of mathematics practitioners as authors, or acknowledgement of their contributions to projects or publications, requires their explicit permission because it implies endorsement of the work. |
| |
| **An ethical mathematics practitioner who is a leader, employer, supervisor, mentor, or instructor follows all of the above items and also...** |
| **IN GENERAL** |
| 34. Maintains a working environment free from intimidation, including discrimination based on personal characteristics; bullying; coercion; unwelcome physical (including sexual) contact; and other forms of harassment. |
| 35. Articulates these ethical responsibilities to mathematics practitioners as well as non-practitioners. |
| 36. Ensures that they enhance, not degrade, the quality of working life. Leaders should consider accessibility, physical safety, psychological well-being, and human dignity of all community members. |
| 37. Does not exploit the offer of a temporary position at an unreasonably low salary and/or an unreasonably heavy workload. |
| **AS A MEMBER OF THE PROFESSION** |
| 38. Recognizes that mathematicians' ethical responsibilities exist and were articulated for the protection and support of the mathematics practitioner, the mathematics user, and the public alike. |
| 39. Encourages and promotes sound and ethical mathematical practice, and exposes incompetent or corrupt mathematical practice. |
| 40. Strives to protect the professional freedom and responsibility of mathematical practitioners who comply with these guidelines. |
| 41. Articulates, applies, and supports policies and processes that reflect the principles of mathematicians' responsibilities. Designing or implementing policies that deliberately or negligently violate, or tend to enable the violation of, mathematicians' responsibilities is ethically unacceptable. |
| 42. Ensures that opportunities are available to mathematics practitioners to help them improve their knowledge and skills in the practice and dissemination of mathematics, in ethical practice, and in their specific fields, and encourages people to take those opportunities. |
| 43. Demonstrates and educates students, employees, and peers on the ethical aspects of their teaching, ethical implications of their work, and the ethical challenges within the practice of mathematics. |
| 44. Takes full responsibility for their contributions to the certification/granting of a degree in mathematics by ensuring the high level and originality of the Ph.D. dissertation work, and sufficient knowledge in the recipient of important branches of mathematics outside the scope of the thesis. |